\documentclass{article}
\usepackage{amsmath,amsfonts,epsfig}
\usepackage{mathptmx}
\newcommand{\bb}{\mathbb}

\newcommand{\cx}{{\bb C}}
\newcommand{\half}{{\bb H}}
\newcommand{\integers}{{\bb Z}}

\newcommand{\reals}{{\bb R}}

\newcommand{\hthree}{{\bb H}^3}
\newcommand{\htwo}{{\bb H}^2}

\newcommand{\tube}{{\mathbb T}}

\renewcommand{\bold}[1]{\medskip \noindent {\bf \boldmath #1
                        }\nopagebreak[4]}

\newcommand{\qed}[1]{\nopagebreak[4]{\tiny \hfill
\fbox{\ref{#1}} \linebreak }\pagebreak[2]}
\newcommand{\bdry}{\partial}

\newcommand{\del}{\partial}
\newcommand{\dirsum}{\oplus}

\newcommand{\chat}{\widehat{\cx}}

\newcommand{\zed}{\integers}

\newcommand{\area}{\operatorname{area}}

\newcommand{\core}{\operatorname{core}}

\newcommand{\Isom}{\operatorname{Isom}}

\newcommand{\PSL}{\operatorname{PSL}}

\newcommand{\sym}{\operatorname{sym}}

\newcommand{\vol}{\operatorname{vol}}

\newtheorem{theorem}{Theorem}[section]
\newtheorem{prop}[theorem]{Proposition}
\newtheorem{lemma}[theorem]{Lemma}

\newcommand{\cC}{{\cal C}}

\newcommand{\cL}{{\cal L}}

\newcommand{\cN}{{\cal N}}

\newcommand{\calM}{{\mathcal M}}

\renewcommand{\hbar}{\bar{{\mathbb H}}^3}

\renewcommand{\core}{C}
\newcommand{\bilip}{\operatorname{bilip}}
\newcommand{\tildeH}{\tilde{{\mathbb H}}^3}

\renewcommand\marginpar[1]{} 

\begin{document}

\title{\bf Geometric inflexibility of hyperbolic cone-manifolds}
\author{Jeffrey Brock\thanks{Research supported by NSF grant
    DMS-1207572.} \  and Kenneth Bromberg\thanks{Research supported by 
    NSF grant DMS-1207873.}}

\date{\today}

\maketitle

\begin{abstract}
We prove 3-dimensional hyperbolic cone-manifolds are geometrically
inflexible: a cone-deformation of a hyperbolic cone-manifold
determines a bi-Lipschitz diffeomorphism between initial and terminal
manifolds in the deformation in the complement of a standard tubular
neighborhood of the cone-locus whose pointwise bi-Lipschitz constant
decays exponentially in the distance from the cone-singularity.
Estimates at points in the thin part are controlled by similar
estimates on the complex lengths of short curves.

\end{abstract}

\section{Introduction}
\label{intro}
In our earlier paper \cite{Brock:Bromberg:inflexible}, we developed an
explicit realization of the qualitative idea that
deformations at infinity of hyperbolic 3-manifolds have effect on the
internal geometry that decays exponentially fast with the depth in the
convex core.  This notion of geometric inflexibility, suggested by
McMullen and exhibited in the restrictive setting of injectivity
bounds, proved sufficiently robust to give a new analytic proof of
Thurston's {\em Double-Limit Theorem} for iteration of pseudo-Anosov mapping
classes and a new ``stand-alone'' proof of the hyperbolization theorem
for 3-manifolds that fiber over the circle with pseudo-Anosov
monodromy.

This paper extends our inflexibility results to the setting where the
change in the geometry is the result of a ``cone-deformation,'' in
which the cone-angle at a closed, geodesic singular locus is changed
while the conformal structure at infinity is held fixed.  Our results
control the best pointwise bi-Lipschitz constant outside of a tubular
neighborhood of the singular locus in domain and range.  The optimal
bi-Lipschitz constant decays to $1$ exponentially fast with the distance from
the tubular neighborhood of the singular locus.

\begin{theorem}\label{coneinflex-intro}
Given $\alpha_0, L, K, \epsilon>0$ and $B>1$ there exists an $R > 0$
and a $d>0$ such that the following holds. Let $(M, g_\alpha)$ be a
geometrically finite hyperbolic cone-manifold with all cone-angles
$\alpha < \alpha_0$ and the length of the singular locus is at most
$L$. Then there exists a one-parameter family of geometrically finite
hyperbolic cone-manifolds $(M, g_t)$ defined for $t \in [0, \alpha]$
so that each component of the singular locus of $(M, g_t)$
has cone-angle $t$ and the conformal boundary
  is the same as the conformal boundary of $(M, g_\alpha)$ so that the
  following holds:
\begin{enumerate}
\item If $U_{\alpha}$ is the $R$-tubular neighborhood of the singular
  locus in $(M, g_\alpha)$ and $U_t$ is a tubular neighborhood of
  the singular locus in $(M, g_t)$ such that $\area(\del U_t) =
  \area(\del U_\alpha)$, then there exists $B$-bi-Lipschitz
  diffeomorphisms 
$$\phi_t: M_\alpha  \backslash U_\alpha \to M_t\backslash U_t$$ 
such that $\phi_t$ is the
  identity map on $M$ in the $\epsilon$-thick part of $M_\alpha$.

\item If $p$
  is in the $\epsilon$-thick part of $(M, g_\alpha)$ then the
  pointwise bi-Lipschitz constant of the maps
$$\phi_t:M_\alpha \to M_t$$
satisfies
$$\log \bilip(\phi_t, p) \leq C_1e^{-C_2 d_\alpha(p, M_\alpha\backslash U_\alpha)}$$
where the constants $C_1$ and $C_2$ depend on the $\alpha_0, L, K,
\epsilon$ and $B$.
\end{enumerate}
\end{theorem}

Similar techniques control the behavior of the complex lengths of
short geodesics in the manifold under the cone deformation, and once
again the distortion decays exponentially in the distance from the
tubular neighborhood of the cone-singularity.
\begin{theorem}\label{conelengthinflex}
Let $M_t = (M, g_t)$ be the one parameter family of geometrically
finite cone-manifolds given by Theorem~\ref{coneinflex-intro}. Let
$\gamma$ be an essential simple closed curve in $M$  and $\gamma_t$
its geodesic representatives in $M_t$. Assume that
$\ell_\alpha(\gamma) < \ell$ for some $\ell>0$. Then there exists
constants $C_1$ and $C_2$ depending on the constants $\alpha_0, L, K,
\epsilon$ and $B$ from Theorem~\ref{coneinflex-intro} and on $\ell$
such that the following holds:
\begin{enumerate}
\item If $\epsilon \leq \ell_\alpha(\gamma) \leq \ell$ then
$$\left| \log \frac{\ell_t(\gamma)}{\ell_\alpha(\gamma)} \right| \leq
C_1 e^{-C_2 d_\alpha(\gamma_\alpha, U_\alpha)}.$$ 

\item If $\ell_\alpha(\gamma) \leq \epsilon/B$ then
$$\left| \log \frac{\ell_t(\gamma)}{\ell_\alpha(\gamma)} \right| \leq
C_1 e^{-C_2 d_\alpha(U^\alpha_{\epsilon}(\gamma), U_\alpha)}.$$ 
\end{enumerate}
\end{theorem}

The idea that complete hyperbolic 3-manifolds are increasingly
inflexible as one takes basepoints deeper and deeper in the convex
core is a natural outgrowth of Mostow and Sullivan rigidity.  McMullen
made this qualitative notion precise in the presence of injectivity
bounds in \cite{McMullen:book:RTM}, but his method made strong use of
geometric limit arguments possible only in the complete setting.  Our
original argument for the complete case in
\cite{Brock:Bromberg:inflexible} shows this pointwise exponential
decay for points outside the thin part, which is an optimal result
(the each tubular thin part is controlled using the complex lengths of
the core geodesics).

Here, the cone-deformation version generalizes the cone-rigidity
theorems of Hodgson-Kerckhoff \cite{Hodgson:Kerckhoff:rigidity} and the second author, and
enhances the bi-Lipschitz metric control away from the cone locus
obtained in \cite{Brock:Bromberg:density} to give explicit decay
estimates in terms of the distance from a standard tubular
neighborhood of the cone locus.

\bold{Inflexibility and ending laminations.}  Geometric inflexibility has provided a range
of new tools to analyze the geometry of and deformation theory of
hyperbolic 3-manifolds.  A key application of the work in the present
paper will be an approach to the geometric classification of finitely
generated Kleinian groups via their {\em ending laminations},
combinatorial invariants that are naturally associated to infinite
volume geometric 'ends' of the convex core of a hyperbolic 3-manifold
with finitely generated fundamental group, which we briefly describe.
The ending lamination records the asymptotics of simple closed curves
on a surface cutting of an end of a hyperbolic 3-manifold, whose
geodesic representatives in the 3-manifold have an {\em a priori}
length bound (and therefore must exit the end of the convex core).

A Theorem of Minsky \cite{Minsky:CKGI} guarantees that for any hyperbolic
3-manifold $M$ in a Bers slice $B_Y$ with the ending lamination
$\lambda$ there is an almost canonical (up to bounded choice at each
stage) sequence of pants decompositions $P_n \to \lambda$ that arises
with uniformly bounded total length $\ell_M(P_n) < L$ in $M$.

The notion of {\em grafting}
\cite{Bromberg:bers,Brock:Bromberg:density} may be employed with a
covering argument similar to that of \cite{Bromberg:Souto:density}, to
allow us to drill the curves in $P_n$ in $M$ with a cone-deformation
that sends the cone angle to zero.  This produces a {\em maximal cusp}
$C_n \in B_Y$, and as the pants decompositions $P_n$ move deeper and
deeper into the convex core, the inflexibility theorem guarantees that
the cone-deformations deform the geometry of $\calM$ in a manner that
decays with the distance of the geodesic representatives of the curves
in $P_n$ from $\calM$.  It follows that $C_n$ limits to $M$, and as
$P_n$ depend only on $\lambda$, the lamination $\lambda$ determines
$M$.  We take up this approach in \cite{BBES:elc}.

\bold{Acknowledgements.}  The authors gratefully acknowledge the
support of the National Science Foundation.

\section{Deformations}
Let $(M, g_t)$ be a one-parameter family of Riemannian manifolds. The time zero derivative of this  family of metrics is given by the formula
$$\eta = \frac{dg_t(v,w)}{dt} \left|_{t=0} \right. = 2g(\eta(v), w).$$
This derivative is a symmetric tensor of type $(1,1)$. We can define a pointwise norm of $\eta$ by fixing an orthonormal basis $\{e_1, \dots, e_n\}$ for $T_p M$ and setting
$$\|\eta\|^2 = \sum_i g(\eta(e_i), \eta(e_i)).$$
As the $L^2$-norm bounds the sup norm we have the inequality
$$\|\eta(v)\| \leq \|\eta\|\|v\|$$
which will be useful in control the change in geometry throughout the flow.

In this paper we will be interested in the case when $(M, g_t)$ is a family of hyperbolic 3-manifolds and the derivative $\eta$ is a {\em harmonic strain field}. Loosely speaking, $\eta$ is harmonic if it locally minimizes the $L^2$-norm. Here is a precise definition. Every point $p$ in $M$ has a chart $U$ and a smooth family of maps $\phi_t: U \to \hthree$ such that on $U$ the hyperbolic metric $g_t$ is the $\phi_t$-pullback of the hyperbolic metric on $\hthree$. For each $q \in U$, $\phi_t(q)$ is a smooth path in $\hthree$ and the time zero tangent vector of this path defines a vector field on $\phi_0(U)$. Let $v$ be the $\phi_0$-pullback of this vector to $U$. If $D$ is the covariant derivative for $g$ then $\eta = \sym Dv$. The infinitesimal change in volume is measured by the trace of $\sym DV$, the {\em divergence} of the vector field. The traceless, symmetric part, $\sym_0 Dv$ is the {\em strain} of $v$ and it measures the infinitesimal change in conformal structure. A vector field is {\em harmonic} if it satisfies the equation
$$D^* D v + 2 v = 0$$
where $D^*$ is the formal adjoint of $D$. While it might be more natural to define $v$ to be harmonic when $D^*Dv=0$ we include the $0$-th order term as we want infinitesimal isometries to be harmonic. This extra term comes from the fact that the Ricci curvature of hyperbolic space is $-2$. We then say that $\eta$ is a {\em harmonic strain field} if $\eta = \sym D v$ where $v$ is a divergence free, harmonic vector field.

On a hyperbolic 3-manifold with boundary, a global bound on the norm of a harmonic strain field leads to exponential decay, in distance from the boundary,  of the pointwise norm in the thick part of the manifold. Before we state the main results from \cite{Brock:Bromberg:inflexible} we make some more definitions. Let $M_t = (M, g_t)$ be a one-parameter family of hyperbolic 3-manifolds. Then $M_t^{\geq \epsilon}$ is the $\epsilon$-thick part of $M_t$, those points where the injectivity radius is $\geq \epsilon$. Here is a key structural theorem from \cite{Brock:Bromberg:inflexible}.
\begin{theorem}
\label{maininflex}
Let $g_t$ be a one-parameter family of hyperbolic metrics on a
3-manifold $M$ with $t \in [a,b]$.  Let $\eta_t$ be the time $t$
derivative of the metrics $g_t$ and let $N_t$ be a family of
submanifolds of $M$ such that $\eta_t$ is a harmonic strain field on
$N_t$. Also assume that
$$\int_{N_t} \|\eta_t\|^2 + \|D_t\eta_t\|^2 \leq K^2$$
for some $K>0$.
Let $p$ be a point in $M$ such that for all $t\in [a,b]$, $p$ is in $M_t^{\geq \epsilon}$ and
$$d_{M_t}(p, M\backslash N_t) \geq	d$$
where $d>\epsilon$.
Then
$$\log \bilip(\Phi_t, p) \leq (t-a)KA(\epsilon)e^{-d}$$
where $\Phi_t$ is the identity map from $M_a$ to $M_t$, 
$$A(\epsilon) = \frac{3e^{\epsilon}\sqrt{2 \vol(B)}}{4\pi f(\epsilon)}$$
and
$$f(\epsilon) = \cosh(\epsilon) \sin(\sqrt 2\epsilon) - \sqrt 2\sinh(R)\cos(\sqrt 2R).$$
\end{theorem}

In the thin part of the manifold, close to a short geodesic, we lack
this level of control. Instead, we control the length the short geodesic
where the change will decay exponentially in the depth of
certain tubular neighborhoods of the short curves. More specifically,
given a short geodesic $\gamma$ we will measure the depth of a tubular
neighborhood $U$ of $\gamma$ where the area of $\del U$ is bounded
below.

\begin{theorem}
\label{mainlengthinflex}
Let $g_t$ be a one-parameter family of hyperbolic metrics on a
3-manifold $M$ with $t \in [a,b]$.  Let $\eta_t$ be the time $t$
derivative of the metrics $g_t$ and let $N_t$ be a family of
submanifolds of $M$ such that $\eta_t$ is a harmonic strain field on
$N_t$. Also assume that
$$\int_{N_t} \|\eta_t\|^2 + \|D_t\eta_t\|^2 \leq K^2$$
for some $K>0$. Let $\gamma_t$ be the geodesics representative on $(M,
g_t)$ of a closed curve $\gamma$ and let $\ell_\gamma(t)$ be the
length of $\gamma$.
\begin{enumerate}
\item Assume that $\gamma_t$ is in $M_t^{\geq \epsilon}$ for all $t\in
  [a,b]$, and that
$$d_{M_t}(\gamma_t, M\backslash N_t) \geq	d.$$
Then
$$\left|\log \frac{\ell_\gamma(b)}{\ell_\gamma(a)}\right| \leq \sqrt{2/3}A(\epsilon)(b-a)Ke^{-d}.$$
\item Assume $\gamma_t$ has a tubular neighborhood $U_t$ of radius
  $\ge R$ and the area of $\del U_t$ is $\ge B$. Also assume that
$$d_{M_t}(U_t, M\backslash N_t) \geq d$$
for all $t\in [a,b]$.
Then
$$\left|\log \frac{\ell_\gamma(b)}{\ell_\gamma(a)}\right| \leq \frac{C(R)(b-a)Ke^{-d}}{\sqrt{B}}$$
where $$1/C(R) =2 \tanh R \left(2 + \frac{1}{\cosh^2 R}\right).$$
\end{enumerate}
\end{theorem}
The Margulis lemma provides an embedded tubular neighborhood about a
sufficiently short geodesic in a hyperbolic 3-manifold: there is a
$\varepsilon$ such that if $\gamma$ is a primitive closed geodesic and
$\operatorname{length}(\gamma)<\epsilon< \varepsilon$ then the
component of the $\epsilon$-thin part that contains $\gamma$ will be a
tubular neighborhood which we denote $U_\epsilon(\gamma)$. This is the
{\em $\epsilon$-Margulis tube} about $\gamma$ and the area of $\del
U_\epsilon(\gamma)$ is bounded below by $\pi \epsilon^2$. In
particular we can apply (2) of the above theorem to such tubes. In
this paper, we will be studying singular hyperbolic manifolds so we
will need to adapt this slightly to find our tubes. 

\section{Cone-manifolds}
We now turn our attention to deformations of hyperbolic
cone-manifolds. We begin with a definition. We let $\tildeH$ be the
set
$$\{(r, \theta, z) | r>0, \theta, z \in \reals\}$$
with the incomplete Riemannian metric 
$$dr^2 + \sinh^2 r d\theta^2 + \cosh^2r dz^2.$$
Then $\tildeH$ is isometric to
the lift to the universal cover of the hyperbolic 
metric on $\hthree \setminus \ell$
where $\ell$ is a complete
geodesic.   For each $\alpha>0$,
let $\hthree_\alpha$ be the metric completion of the quotient of
$\tildeH$ under the isometry $(r, \theta, z) \mapsto (r,
\theta+\alpha,z)$. Note that $\hthree_\alpha$ is a topological
ball. Let $N$ be a compact 3-manifold with boundary and $g$ a complete
metric on the interior of $N$. The metric $g$ is a {\em hyperbolic
  cone-metric} if every point in the interior of $N$ has a
neighborhood isometric to a neighborhood of a point in
$\hthree_\alpha$ for some $\alpha >0$. The pair $(N, h)$ is a {\em
  hyperbolic cone-manifold}. Let $\cC$ be the subset of $N$ where the
metric $h$ is singular.  Then $\cC$ will be a collection of isolated
simple curves in $N$. In this paper we will assume that $\cC$ is
compact which implies that it is a finite collection of disjoint
simple closed curves.

Let $c$ be a component of $\cC$. Then there is a unique $\alpha>0$ such that each point $p$ in $c$ has a neighborhood isometric to the neighborhood of a singular point in $\hthree_\alpha$. This $\alpha$ is the cone-angle of the component $c$.

Recall that $\hthree$ is naturally compactified by $\chat$. The union
is a closed 3-ball and isometries of $\hthree$ extend continuously to
conformal automorphisms of $\chat$. Let $\del_0 N$ be the components
of $\del N$ that are not tori. Then $(N, g)$ is a {\em geometrically
  finite cone-manifold} if each point $p$ in $\del_0 N$ has a
neighborhood $V$ in $N$ and a chart $\phi: V \to \hbar$ such that
$\phi$ restricted to $V \cap \operatorname{int}(N)$ is an isometry and
$\phi$ restricted to $V \cap \del N$ is a map into $\del \hbar =
\chat$. Note that the restriction of the charts to $\del_0 N$ defines
an atlas for a conformal structure on $\del_0 N$. In fact, as we will be important in the next section, this conformal atlas determines a {\em complex
projective structure} on $\del_0 N$.

\begin{theorem}
\label{conedeformation}
Given $\alpha_0, L, K, \epsilon>0$ and $B>1$ there exists an $R > 0$
and a $d>0$ such that the following holds. Let $M_\alpha= (M, g_\alpha)$ be a
geometrically finite hyperbolic cone-manifold with all cone-angles
$\alpha < \alpha_0$, each component of the singular locus has an embedded tubular neighborhood of radius $R$ and the length of the singular locus is at most
$L$. Then there exists a one-parameter family of geometrically finite
hyperbolic cone-manifolds $M_t = (M, g_t)$ defined for $t \in [0, \alpha]$
with the following properties:
\begin{enumerate}
\item Each component of the singular locus of $M_t$
has cone-angle $t$ and the conformal boundary
  is the same as the conformal boundary of $M_\alpha$.

\item The derivative $\eta_t$ of $g_t$ is a family of harmonic strain
  fields outside of a radius $\sinh^{-1} (1/\sqrt{2})$
  neighborhood of the singular locus.

\item Let $U_{\alpha}$ be the $R$-tubular neighborhood of the singular
  locus in $M_\alpha$ and let $U_t$ be a tubular neighborhood of
  the singular locus in $M_t$ such that $\area(\del U_t) =
  \area(\del U_\alpha)$. Then
$$\int_{M_t \backslash U_t} \|\eta_t\|^2 + \|D_t \eta_t\|^2 \leq K^2.$$

\item There exists $B$-bi-Lipschitz diffeomorphisms $\phi_t: M_\alpha
  \backslash U_\alpha \to M_t\backslash U_t$ such that $\phi_t$ is the
  identity map on $M$ in the $\epsilon$-thick part of $M_\alpha$.

\item If $p \in (M_\alpha\backslash U_\alpha)^{\ge \epsilon}$ then $p \in (M_t\backslash U_t)^{\ge \epsilon/B}$ and $d_t(p, U_t) \geq d_\alpha(p,
  U_\alpha)/B$.

\item If $\gamma$ is a closed curve in $M$ then $d_t(\gamma_t, U_t)
  \geq d_\alpha(\gamma_\alpha, U_\alpha)/B -d$.

\item If $\gamma$ is a closed curve in $M$ with $\ell_\alpha(\gamma) < \epsilon/B$ then
$$d_t(U^t_{\epsilon}(\gamma), U_t) \geq \frac{d_\alpha(U^\alpha_{\epsilon}(\gamma), U_\alpha)}{B}-d.$$
\end{enumerate}
\end{theorem}

\bold{Proof.} Statements (1)-(4) are proven in \cite{Bromberg:long}
(see Theorem 5.3 and its proof). When the singular locus is
sufficiently short this was proven in \cite{Bromberg:Schwarzian,
  Brock:Bromberg:density} building on Hodgson and Kerckhoff's
foundational work on deformations of hyperbolic cone-manifolds in
\cite{Hodgson:Kerckhoff:rigidity, Hodgson:Kerckhoff:bounds,
  Hodgson:Kerckhoff:shape}.

Statement (5) follows directly from (4). Statements (6) and (7) are
more difficult. To prove them we need to modify the
metrics $g_\alpha$ and $g_t$ in $U_\alpha$ and $U_t$ so that they are
complete metrics of pinched negative curvature and by then extending the map
$\phi_t$ to a bi-Lipschitz map for these new metrics.

The construction of such metrics is straightforward: they are doubly
warped products using cylindrical coordinates.
Given an $r_0>0$ define a metric on $\reals^3$ by
$$dr^2 + f_{r_0}(r)^2 d\theta^2 + g_{r_0}(r)^2 dz^2$$
where $f_{r_0}(r)$ and $g_{r_0}(r)$ are convex functions with
$f_{r_0}(r) = \sinh r$ and $g_{r_0}(r) = \cosh r$ for $r \in [{r_0}/2,
{r_0}]$ and $f_{r_0}(r) = g_{r_0}(r) = \frac12 e^r$ for $r \leq
{r_0}/4$. We can also assume that $\sinh r \leq f_{r_0}(t) \leq
\frac12 e^r$ and $\frac12 e^r \leq g_{r_0}(r) \leq \cosh r$.  When
$r\geq {r_0}/2$ or $r \leq {r_0}/4$ then this metric is
hyperbolic. For $r \in ({r_0}/4, {r_0}/2)$ the sectional curvature
will be pinched within $\delta$ of $-1$ where $\epsilon$ only
depends on ${r_0}$ and $\delta \to 0$ as ${r_0} \to \infty$. Details
of this calculation can be found in Section 1.2 of \cite{Kojima:cone}
where the construction is attributed to Kerckhoff.

The map $(r,\theta,z) \mapsto (r, \theta + x, z + y)$ is an
isometry in this metric. If we take the quotient of the set of points
with $r \in (-\infty,r_0]$ by isometries $(r, \theta, z) \mapsto (r,
\theta +t, z)$ and $(r, \theta + x, z + \ell)$ we get a complete
metric on $T^2 \times (-\infty, r_0]$. If $r_0=R_t$ is the tube radius
of $U_t$ and $\ell + \imath x$ is the complex length of the singular
locus of $(M, g_t)$ then the $R_t/2$-neighborhood of the boundary is
isometric to the $R_t/2$-neighborhood of $\del U_t$. We then define
$g'_t$ on $U_t$ by replacing the original metric with the above
metric. Since the two metrics agree in a collar neighborhood of $\del
U_t$ the metric $g'_t$ is smooth and $g'_t$ is a complete metric on
$M$ with sectional curvature within $\delta$ of $-1$.

We now construct a bi-Lipschitz diffeomorphism $\phi'_t:(M, g'_\alpha)
\to (M, g'_t)$ by extending the map $\phi_t$ from (4). The original
map $\phi_t$ restricted to $\del U_\alpha$ is a $B$-bi-Lipschitz
diffeomorphism from $\del U_\alpha$ to $\del U_t$. This map can then
be extended to a map on $(U_\alpha, g'_\alpha)$ in the obvious
way. Namely there are nearest point projections of $(U_\alpha,
g'_\alpha)$ and $(U_t, g'_t)$ onto $\del U_\alpha$ and $\del U_t$
respectively. Then on $U_\alpha$, $\phi'_t$ is the unique map that
commutes with these projections and that takes a point distance $r$
from $\del U_\alpha$ to a point distance $r$ from $\del U_t$. We need
to calculate the bi-Lipschitz constant of this map.

To do so we make a few observations. First the functions $f_R(r)$ and
$g_R(r)$ converge uniformly to $\frac12 e^r$ as $R \to \infty$. Second
we note that by construction the derivative of the map is an isometry
in the $r$-direction. For a vector $v$ tangent to the tori of fixed
$r$-coordinate a direction calculation shows that
$$\frac{1}{B} \frac{f_{R_t}(r')}{f_{R_t}(R_t)} \frac{g_{R_\alpha}(R_\alpha)}{g_{R_\alpha}(r)} \|v\| \leq \|\left(\phi'_t\right)_* v\| \leq B \frac{g_{R_t}(r')}{g_{R_t}(R_t)} \frac{f_{R_\alpha}(R_\alpha)}{f_{R_\alpha}(r)} \|v\|$$
where $R_\alpha -r = R_t - r'$. Therefore the map is $B'$-bi-Lipschitz
where $B'$ is the maximum of the factor on the right side of the
inequality and the inverse of the factor on left side of the
inequality. Since the functions $f_R(r)$ and $g_R(r)$ converge
uniformly to $\frac12 e^r$, the quotients $f_R(r_1)/f_R(r_0)$ and
$g_R(r_1)/g_R(r_0)$ converge uniformly to $e^{r_1 - r_0}$.  By Theorem
2.7 of \cite{Hodgson:Kerckhoff:bounds} the length of the singular
locus is an increasing function of $t$. This implies that $R_t$ is a
decreasing function in $t$ and therefore the bi-Lipschitz constant,
$B'$, depends only on $B$ and $R$.

By the Morse Lemma (see e.g. \cite{Bridson:Haefliger:npc}) the $\phi_t$-image of a geodesic is contained in the $d$-neighborhood of a geodesic where $d$ only depends on $B'$ and the curvature bounds of the modified metric (which we have uniformly controlled). In particular, $\phi_t(\gamma_\alpha)$ is contained in the $d$-neighborhood of $\gamma_t$. Since $\phi_t$ is $B$-bi-Lipschitz on $M_\alpha \backslash U_\alpha$ and $\phi_t(U_\alpha) = U_t$ we have $d_t(\phi_t(\gamma_\alpha), U_t) \geq d_\alpha(\gamma_\alpha, U_\alpha)/B$ and therefore $d_t(\gamma_t, U_t) \geq d_\alpha(\gamma_\alpha, U_\alpha)/B-d$ which is (6).

For (7) we choose $\epsilon$ such that the $B\epsilon$ is less than than Margulis constant for manifolds with curvature pinched between $-1 - \delta$ and $-1 + \delta$. Then if $\ell_\alpha(\gamma) < \epsilon/B$ we have that $\ell_t(\gamma) < \epsilon < B\epsilon$ and both $U^t_{B\epsilon}(\gamma)$ and $U^t_{\epsilon/B}(\gamma)$ will be embedded tubular neighborhoods. Furthermore we have $U^t_{\epsilon/B}(\gamma) \subseteq \phi_t(U^\alpha_\epsilon(\gamma)) \subseteq U^t_{B\epsilon}(\gamma)$. By \cite{Brooks:Matelski:collars} the width of the collar $U^t_{B\epsilon}(\gamma) - U^t_{\epsilon/B}(\gamma)$ is bounded by a constant that is independent of $\ell_t(\gamma)$. This gives uniform control of the distance between $\phi_t(U^\alpha_\epsilon(\gamma))$ and $U^t_\epsilon(\gamma)$ and then (7) follows in a similar manner as (6).
\qed{conedeformation}

We can now prove the bi-Lipschitz inflexibility theorem for
cone-manifolds.
\begin{theorem}\label{coneinflex}
  Let $M_t = (M,g_t)$ be the one-parameter family of geometrically
  finite cone-manifolds given by Theorem~\ref{conedeformation}. If $p$
  is in the $\epsilon$-thick part of $(M, g_\alpha)$ then the
  pointwise bi-Lipschitz constant of the maps
$$\phi_t:M_\alpha \to M_t$$
satisfies
$$\log \bilip(\phi_t, p) \leq C_1e^{-C_2 d_\alpha(p, U_\alpha)}$$
where the constants $C_1$ and $C_2$ depend on the $\alpha_0, L, K,
\epsilon$ and $B$ as in Theorem~\ref{conedeformation}.
\end{theorem}

\bold{Proof.} We apply Theorem \ref{maininflex} to $M_t$ with $N_t = M_t\backslash U_t$. By (2) of Theorem \ref{conedeformation} the derivative $\eta_t$ of $M_t$ is a harmonic strain field on $N_t$ and by (3) we have that
$$\int_{N_t} \|\eta_t\|^2 +\|D_t \eta_t \|^2 \leq K^2.$$
Let $B>1$ be the bi-Lipschitz constant given by (4) and then, by (5), a point $p \in M^{\ge \epsilon}_\alpha$ will be in $M^{\ge \epsilon/B}_t$ and
$$d_t(p, U_t) \ge d_\alpha(p, U_\alpha)/B.$$
The result then follows from Theorem \ref{maininflex} with $C_1 = \alpha KA(\epsilon/B)$ and $C_2 = 1/B$.
\qed{coneinflex}

Next we state and prove the length inflexibility statement.
\begin{theorem}\label{conelengthinflex}
Let $M_t = (M, g_t)$ be the one parameter family of geometrically finite cone-manifolds given by Theorem~\ref{conedeformation}. Let $\gamma$ be an essential simple closed curve in $M$  and $\gamma_t$ its geodesic representative in $M_t$. Assume that $\ell_\alpha(\gamma) < \ell$ for some $\ell>0$. Then there exists constants $C_1$ and $C_2$ depending on the constants $\alpha_0, L, K, \epsilon$ and $B$ from Theorem~\ref{conedeformation} and on $\ell$ such that the following holds.
\begin{enumerate}
\item If $\epsilon \leq \ell_\alpha(\gamma) \leq \ell$ then
$$\left| \log \frac{\ell_t(\gamma)}{\ell_\alpha(\gamma)} \right| \leq C_1 e^{-C_2 d_\alpha(\gamma_\alpha, U_\alpha)}.$$

\item If $\ell_\alpha(\gamma) \leq \epsilon/B$ then
$$\left| \log \frac{\ell_t(\gamma)}{\ell_\alpha(\gamma)} \right| \leq C_1 e^{-C_2 d_\alpha(U^\alpha_{\epsilon}(\gamma), U_\alpha)}.$$
\end{enumerate}
\end{theorem}

\bold{Proof.} As in the proof of Theorem \ref{coneinflex} we let $N_t = M_t/U_t$ and then by (2) and (3) of Theorem \ref{conedeformation} the derivative of $M_t$ on $N_t$ is a harmonic strain field $\eta_t$ with
$$\int_{N_t} \|\eta_t\|^2 + \|D_t \eta_t\|^2 \le K^2.$$
If $B>1$ is the bi-Lipschitz constant from (4) then by (6) there is a constant $d>0$ such that
$$d_t(\gamma_t, U_t) \ge d_\alpha(\gamma_\alpha, U_\alpha)/B - d.$$
The first inequality the follows from (1) of Theorem \ref{mainlengthinflex} with $C_1 = \sqrt{2/3}A(\epsilon/B)\alpha K e^{-d}$ and $C_2 = 1/B$.

The second inequality is proved similarly but we use (7) of Theorem \ref{conedeformation} instead of (6).
\qed{conelengthinflex}

\section{Schwarzian derivatives}
As was noted when defining geometrically finite hyperbolic cone-manifolds, the conformal boundary of a hyperbolic cone-manifold also has a projective structure. While the conformal boundary will be fixed throughout the deformations given by Theorem \ref{conedeformation}, the projective structure will vary. The variation in a projective structure is measured by the Schwarzian derivative and in this section we will use our inflexibility theorems to control the size of the Schwarzian derivative.

We very briefly discuss projective structures and the Schwarzian
derivatives. For more detail see Section 6 of
\cite{Brock:Bromberg:inflexible}.  A projective structure on a surface
is  $(G,X)$-structure where $G = PSL_2\cx$ and $X = \chat$. In a
projective chart the derivative of a smooth 1-parameter family of
projective structures is a conformal vector field. Using the chart
this is a vector field $v$ on a domain in $\chat$. At each point there
is a unique projective vector field that best approximates $v$. In
such a way $v$ defines a map from the domain in $\chat$ to $sl_2\cx$
the Lie algebra of projective vector fields. The derivative of this
map is the {\em Schwarzian derivative} of the deformation and it
naturally identified with a holomorphic quadratic differential on the
conformal structure. 

Given two projective structures we define the notion of a {\em projective map} between them in the usual way via charts. For example a round disk in $\chat$ inherits a projective structure as a subspace of $\chat$.  On a arbitrary projective structure $\Sigma$ a {\em round disk} is a projective map from a round disk to $\Sigma$. Note that we don't assume that this map is an embedding. Every round disk in $\chat$ bounds a half space $\hthree$. If $\Sigma$ is the projective boundary of a hyperbolic 3-manifold $M$ then a round disk in $\Sigma$ bounds a half space in $M$ if there is an isometry from a half space in $\hthree$ into $M$ that extends to a projective map on the boundary round disk. We will need the following lemma about round disks.
\begin{lemma}\label{embedded}
Let $\Sigma$ be projective structure with trivial holonomy. The every round disk is embedded.
\end{lemma}

\bold{Proof.} Let $\tilde\Sigma$ be the universal cover of $\Sigma$. Recall that there is a projective developing map $D:\tilde\Sigma \to \chat$ and a representation $\rho: \pi_1(\Sigma) \to PSL_2\cx$ such that $D \circ \gamma = \rho(\gamma) \circ D$ where the action of $\gamma$ in the left side of the inequality is by deck transformations. By assumption the holonomy representation $\rho$ is the trivial representation. 

Let $U$ be a round disk in $\chat$ and $\phi:U \to \Sigma$ projective map. Let $\tilde\phi: U \to \tilde\Sigma$ be the lift of $\phi$. Then $D \circ \tilde\phi$ is a projective map of $U \subset \chat$ into $\chat$. Since $D \circ \tilde\phi$ is the restriction of an element of $PSL_2\cx$ it is an embedding and hence $\tilde\phi$ is an embedding. If $\phi$ is not an embedding then there exists $x, y \in U$ such that $\phi(x) = \phi(y)$. Since $\tilde\phi(x) \neq \tilde\phi(y)$ there must be a $\gamma \in \pi_1(\Sigma)$ such that $\gamma(\tilde\phi(x)) = \tilde\phi(y)$. Since $D \circ \gamma(\tilde\phi(y)) = \rho(\gamma) \circ D(\tilde\phi(y))$ and $\rho(\gamma)$ is the identity we have $D(\tilde\phi(x)) = D(\tilde\phi(y))$. Since $D \circ \tilde\phi$ is injective this is a contradiction and hence $\phi$ is injective.
\qed{embedded}

We now state the main inflexibility theorem for Schwarzian derivatives from \cite{Brock:Bromberg:inflexible}. As the projective structure is at infinity we can't measure its distance from the cone singularity. Instead we assume that each round disk in the projective structure bounds a half space in the manifold and then measure the distance to the half space. 
\begin{theorem}\label{mainschwarzinflex}
Let $g_t$, $t \in [a,b]$, be a one-parameter family
  of hyperbolic metrics on the interior of a 3-manifold $M$ with
  boundary.  Let $\eta_t$ be the time $t$ derivative
  of the metrics $g_t$ and let $N_t$ be a family of submanifolds of
  $M$ with compact boundary such that $\eta_t$ is a harmonic strain
  field on $N_t$. Also assume that
$$\int_{N_t} \|\eta_t\|^2 + \|D_t \eta_t\|^2 \leq K^2$$
for some $K>0$.  Let $S$ be a component of $\del M$ such that each
hyperbolic metric $g_t$ extends to a fixed conformal structure $X$ on
$S$ and a family of projective structures $\Sigma_t$ on $S$. Assume
that at every embedded round disk in $\Sigma_t$ bounds an embedded
half space $H$ in $N_t$ and that
$$d_{M_t}(H, M \backslash N_t) \geq d$$
for some $d>0$.  Then
$$d(\Sigma_a, \Sigma_b) \leq CKe^{-d}$$
where $C$ is a constant depending on the sup-norm of the Schwarzian derivative of the quadratic differential from the unique Fuchsian projective structure with conformal structure $X$ and the
injectivity radius of $X$.
\end{theorem}

To apply this result we need to know that round disks in the projective boundary of a hyperbolic cone-manifold bound a half spaces.
\begin{lemma}\label{rounddisks}
  Let $M$ be the non-singular part of a 3-dimensional hyperbolic
  cone-manifold. Then every round disk on the projective boundary of
  $M$ extends to a half-space in $M$, and if the disk is embedded the
  half space is embedded.
\end{lemma}

\bold{Proof.} In Lemma 3.3 of \cite{Bromberg:Schwarzian} it is shown
that every embedded round disk extends to an embedded half space so we
only need to show that every (possibly immersed) round disk extends to a half space. To do this we would like to apply the lemma to the universal cover $\tilde M$ of the non-singular part of the hyperbolic cone-manifold. We first observe that if $\tilde\Sigma$ is a component of the projective boundary then its holonomy representation will be trivial so by Lemma \ref{embedded} every round disk in $\tilde\Sigma$ will be embedded. On the other hand, $\tilde M$ is not a hyperbolic cone-manifold in the sense that is used in the proof of Lemma 3.3 of \cite{Bromberg:Schwarzian} so we will briefly review the proof to see that it applies in our situation.

A hyperbolic half space $H \subset \hthree$ is foliated by constant curvature planes $P_d$ where $P_d$ is the locus of points distance $d$ from the hyperbolic plane that bounds $H$. Let $H_d \subset H$ be the union of the $P_t$ with $t > d$. Let $U$ be a round disk in $\tilde\Sigma$ whose closure is compact. Using a compactness argument we can extend the round disk to $H_d$ for some for some large $d$. We identify $H_d$ with its image in $\tilde M$. When $d>0$ the boundary of $H_d$ is strictly concave so $\tilde M\backslash H_d$ is strictly convex and therefore the closure of $H_d$ is embedded in $\tilde M$ if $d>0$. This implies that we can extend the round disk to $H_0$. The hyperbolic plane that 
\qed{rounddisks}

If $\Sigma$ is the projective boundary of a hyperbolic cone-manifold $M$ we define its neighborhood $\cN(\Sigma)$ to be the union of all half-spaces that are
bounded by round disks in $\Sigma$. Since two half-spaces in $M$ will
intersect if and only if their boundary round disks intersect,
disjoint components of the projective boundary will determine disjoint
neighborhoods.

Thurston parameterized the space of projective structures on a surface
$S$ by the product of the Teichm\"uller space and the space of
measured laminations. In his proof he extends a projective structure
to a hyperbolic structure on $\Sigma \times [0,\infty)$ where the
boundary is a locally concave pleated surface (or a locally convex
pleated surface if it is embedded in a larger
manifold). Lemma~\ref{rounddisks} essentially shows that this
hyperbolic structure constructed by Thurston is our neighborhood
$\cN(\Sigma)$. We now state Thurston's result in a form that will be
useful to us. For a proof see \cite{Kamishima:Tan}.

\begin{theorem}[Thurston]\label{thurston}
  Each neighborhood $\cN(\Sigma)$ is homeomorphic to $\Sigma \times
  (0, \infty)$. If the singular locus doesn't intersect the boundary
  of $\cN(\Sigma)$ then the boundary is a locally convex pleated
  surface.
\end{theorem}

Our inflexibility theorems will be vacuous if the singular locus is on
the boundary of $\cN(\Sigma)$ so we can effectively assume that this is not the case and that the
boundary of $\cN(\Sigma)$ is a locally convex pleated surface.

The {\em convex core} of a complete manifold of pinched negative curvature is the smallest convex subset whose inclusion is a homotopy equivalence. As the non-singular part of a cone-manifold is not complete we need to be more careful in how we define the convex core. The following lemma will be essential.

\begin{lemma}\label{coneconvexcore}
  Let $(M, g)$ be the non-singular part of a 3-dimensional hyperbolic
  cone-manifold and let $(M, g')$ be a complete Riemannian metric on
  $M$ with pinched negative curvature such that $g=g'$ on
  $\cN(\Sigma)$. Then $M \backslash \cN(\Sigma)$ is the convex core of
  $(M, g')$.
\end{lemma}

\bold{Proof.} By Theorem~\ref{thurston} the manifold $M$ deformation
retracts onto $M \backslash \cN(\Sigma)$ so the inclusion of $M
\backslash \cN(\Sigma)$ into $M$ will be a homotopy equivalence. The
boundary of $M \backslash \cN(\Sigma)$ will be locally convex in
$(M,g)$ and therefore also in $(M, g')$.  This implies that
$M\backslash \cN(\Sigma)$ is a convex sub-manifold in $(M,g')$ whose
inclusion is a homotopy equivalence and therefore the convex core is
contained in $M \backslash \cN(\Sigma)$.

Next we show that the pleating locus of the pleated surfaces bounding
$M \backslash \cN(\Sigma)$ must be contained in the convex core. To
see this we first note that any closed geodesic is in the convex
core. The pleating locus can be approximated by closed geodesics so it
must also be in the convex core.  

Finally the join of anything in the convex core will also be in the
convex core. Since the join of the pleating locus will contain the
pleated surface we have that  $\bdry (M \backslash
\cN(\Sigma))$ lies in the convex core so $M\backslash
\cN(\Sigma)$ lies in the convex core. \qed{coneconvexcore}

Given this lemma, it is natural to define the convex core of a
hyperbolic cone-manifold by $C(M) = M \backslash \cN(\Sigma)$.
For this definition to be useful we need to know that the image of the convex core under a bi-Lipschitz map will be uniformly close in the Hausdorff metric to the convex core of the image manifold. This will follow from the following proposition which is due to McMullen when the manifold is hyperbolic. The general case requires work of Anderson and Bowditch.
\begin{prop}\label{coresclose}
  Given $B>1$ and $\epsilon \in (0,1)$ there exists $d>0$ such that
  the following holds. Let $g_0$ and $g_1$ be complete Riemannian
  metrics on a manifold $M$ with sectional curvatures in
  $(-1-\epsilon,-1 + \epsilon)$ and let $\phi: (M, g_0)
  \to (M, g_1)$ be $B$-bi-Lipschitz. Then then Hausdorff
  distance between $C(M,g_1)$ and $\phi(C(M,g_0))$ is less than $d$.
\end{prop}

The final piece we need to prove our Schwarzian inflexibility theorem is a version of the deformation theorem for cone-manifolds that controls the distance from the standard neighborhood of the singular locus to the convex core boundary. It will be convenient to restate part of the original deformation theorem, Theorem~\ref{conedeformation}.

\begin{theorem}
\label{conedeformationwithcore}
Given $\alpha_0, L, K>0$ and $B>1$ there exists an $R > 0$ such that
the following holds. Let $(M, g_\alpha)$ be a geometrically finite
hyperbolic cone-manifold with all cone-angles $\alpha < \alpha_0$ and
with singular locus of length at most $L$. Then there exists a one-parameter family of geometrically finite hyperbolic cone-manifolds $(M, g_t)$ defined for $t \in [0, \alpha]$ with the following properties:
\begin{enumerate}
\item All cone angles of $(M, g_t)$ are $t$ and the conformal boundary
  is the same as the conformal boundary of $(M, g_\alpha)$.

\item The derivative $\eta_t$ of $g_t$ is a family of harmonic strain
  fields outside of a radius $\sinh^{-1} 1/\sqrt{2}$ neighborhood of
  the singular locus. 

\item Let $U_{\alpha}$ be the $R$-tubular neighborhood of the singular
  locus in $(M, g_\alpha)$ and let $U_t$ be a tubular neighborhood of
  the singular locus in $(M, g_t)$ such that $\area(\del U_t) =
  \area(\del U_\alpha)$. Then
$$\int_{M_t \backslash U_t} \|\eta_t\|^2 + \|D_t \eta_t\|^2 \leq K.$$

\item \label{thickcore} Let $X$ be a component of the conformal boundary and $\Sigma_t$ the projective structure on $X$ induced by $(M,g_t)$. Then
 $$d(U_t, \cN(\Sigma_t)) \geq d(U_\alpha, \cN(\Sigma_\alpha))/B - d.$$

\end{enumerate}
\end{theorem}

\bold{Proof.} Except for \eqref{thickcore} this is exactly the same Theorem~\ref{conedeformation}. To prove \eqref{thickcore} we would like to apply Proposition~\ref{coresclose} but since our metrics are incomplete we cannot do so directly. We will use the same trick that we used in the proof of Theorem~\ref{conedeformation} and replace the metrics $g_\alpha$ and $g_t$ with complete metrics of pinched negative curvature, $g'_\alpha$ and $g'_t$ and then use the extended $B$-bi-Lipschitz diffeomorphism $\phi'_t$ from $(M, g'_\alpha)$ to $(M, g'_t)$. We then apply Proposition~\ref{coresclose} which shows that
$$Bd(U_t, M\backslash C(M, g'_t)) + d \geq d(U_\alpha, M\backslash C(M, g'_\alpha)).$$
Note that we can assume that $U_\alpha$ is contained in $C(M, g_\alpha)$ for otherwise \eqref{thickcore} is vacuous. The inequality then follows from Lemma~\ref{coneconvexcore}. \qed{conedeformationwithcore}

We can now apply Theorems~\ref{mainschwarzinflex} and~\ref{conedeformationwithcore} to get our Schwarzian inflexibility theorem for cone-manifolds.
\begin{theorem}\label{coneschwarzinflex}
  Given $\alpha_0, L, K>0$ and $B>1$ there exists an $R > 0$ such that
  the following holds. Let $(M, g_\alpha)$ be a geometrically finite
  hyperbolic cone-manifold with all cone-angles $\alpha < \alpha_0$,
  singular locus of length at most $L$ and tube radius of the singular
  locus at least $R$. Let $M_t = (M,g_t)$ be the one-parameter family
  of geometrically finite cone-manifolds given by
  Theorem~\ref{conedeformationwithcore}. Let $\Sigma_t$ be a component
  of the projective boundary of the $M_t$ with underlying conformal
  structure $X$.  Then
$$d(\Sigma_\alpha, \Sigma_t) \leq CK e^{-d(U_\alpha,  \cN(\Sigma_\alpha))/B - d}$$
where $U_\alpha$ is the tubular neighborhood of the singular locus of radius $R_0$ and $C$ is a constant depending on $\|\Sigma_\alpha\|_F$ and the injectivity radius of $X$.
\end{theorem}

\bold{Proof.}
We apply Theorem \ref{mainschwarzinflex} to $M_t$ where the convex cores $C(M_t)$ play the role of the submanifolds $N_t$.  Every half space $H$ bounding a round disk in $\Sigma_t$ will be contained in $\cN(\Sigma_t)$ so by (4) of Theorem \ref{conedeformationwithcore} there exists $d>0$ such that
$$d(U_t, \cN(\Sigma_t)) \ge d(U_\alpha, \cN(\Sigma_\alpha))/B - d.$$
The theorem then follows from (3) of Theorem \ref{conedeformationwithcore} and Theorem \ref{mainschwarzinflex}.
\qed{coneschwarzinflex}

\bibliographystyle{math}
\bibliography{math}

{\sc \small
 \bigskip

\noindent Brown University \bigskip

\noindent University of  Utah
}

\end{document}